\documentclass[11pt,a4paper,onecolumn]{article}

\usepackage[utf8x]{inputenc} 
\usepackage{mathrsfs}
\usepackage{amsfonts}
\usepackage{amsmath}
\usepackage[caption=false,font=footnotesize]{subfig}
\usepackage{amssymb}

\usepackage{ucs}
\usepackage{graphicx}

\usepackage{keyval}
\usepackage{enumerate}
\usepackage{caption}

\usepackage{mathptmx}
\usepackage{color}
\usepackage[normalem]{ulem}
\usepackage[colorlinks=true]{hyperref}
\usepackage{soul}
\usepackage{breakcites}

\hypersetup{
pdftitle = Model-Free Based Digital Control for Magnetic Measurements,
pdfauthor = Loïc Michel - Ghibaudo - Messal - Kedous Lebouc - Boudinet - Blache - Labonne,
pdfsubject = ArXiv paper,
pdfkeywords = digital control - model-free control - magnetic hysteresis - magnetic measurements,
colorlinks = false,       
} 

\begin{document}


\begin{center}
	\begin{huge}
\vspace{0.5cm}
\textbf{Model-Free Based Digital Control \\ for Magnetic Measurements}\\
	\end{huge}
\vspace{0.5cm}
L. Michel, O. Ghibaudo, O. Messal, A. Kedous-Lebouc,\\ C. Boudinet, F. Blache, A. Labonne. \\
\vspace{0.2cm}
Univ. Grenoble Alpes, G2Elab, F-38000 Grenoble, France\\CNRS, G2Elab, F-38000 Grenoble, France.\\	
\end{center}

\noindent \textbf{Abstract -} This paper presents a novel digital control strategy successfully implemented for a soft magnetic material characterization bench (Epstein frame type). The main objective is to control the magnetic induction waveform whatever the applied excitation and the material under study. Given the nonlinear nature of the magnetization  curves of magnetic materials, an original model-free based control technique is considered. Special mention should be made of the interesting dynamic properties in closed-loop against the changes of the operating point related basically to the hysteresis form. The operation and the performances of the digital control method are illustrated in different working conditions through both simulation and experimental measurements.\newline

\noindent \textbf{Keyword -} Digital control, model-free control, magnetic hysteresis, magnetic measurements.

%

\section{Introduction}
\begin{figure}[!htb]
	\centering
		\includegraphics[width=14.5cm]{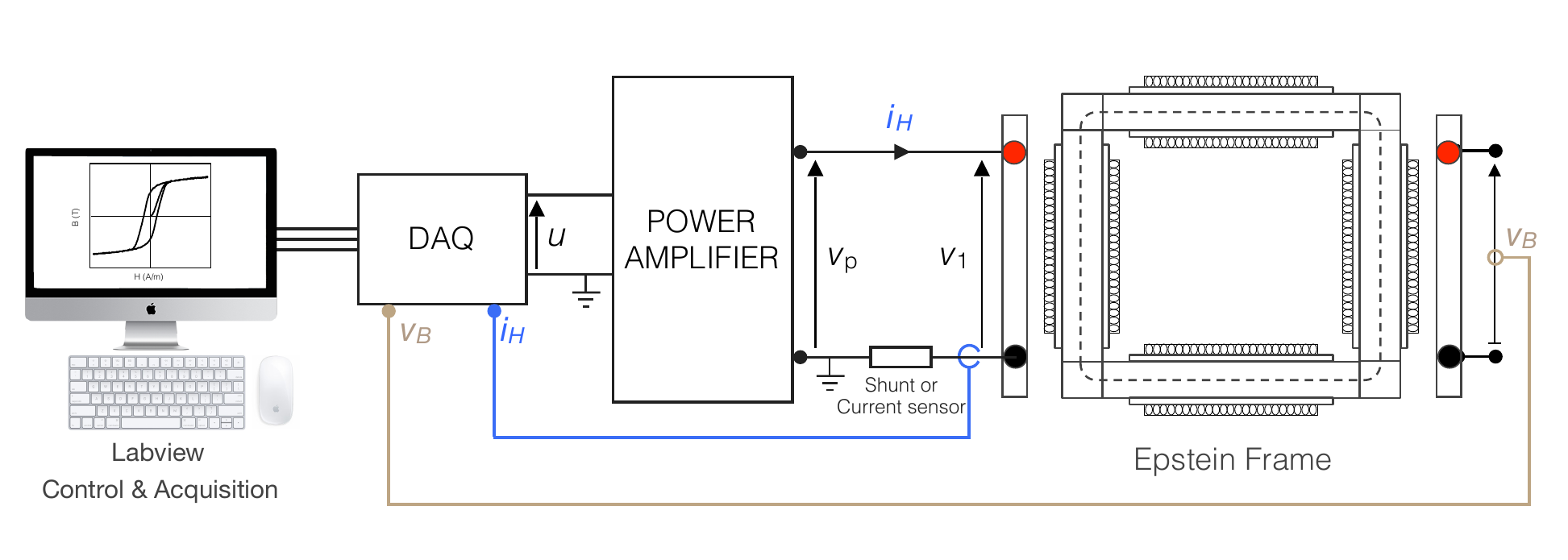}
		\caption{Scheme of the magnetic characterization bench.}
		\label{CaracBench}
\end{figure}

To study of the ferromagnetic materials behavior, namely in their practical working conditions in their main applications in electrical engineering, and for certain modeling requirements, we are often led to do experimental measurements under different induction waveforms (e.g., sine, triangular,...). The experiments are usually performed by means of the standardized Epstein frame device. Its operating principle consists of creating a magnetic field $H$ within the material (according to the Ampère's law) thanks to an excitation current $i_H$ flowing through the Epstein frame primary coil. The material accordingly provides a magnetic response to the applied excitation which physically corresponds to the magnetic induction $B$ obtained after integrating the measured induced voltage $v_B$ in the secondary coil. \\
In practice, the excitation coil is supplied with a linear power supply and indirectly imposes the secondary voltage $v_B$, image of the induction (see Fig. \ref{CaracBench}).\\

\begin{figure}[!htb]
	\centering
		\includegraphics[width=8.8cm]{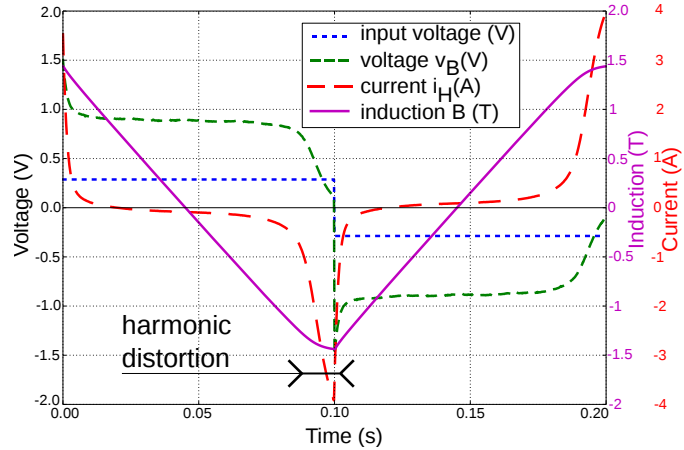}
		\caption{Epstein frame measurements made on a non-oriented SiFe steel under 5 Hz quasi-triangular waveform without feedback. The voltage output reference (input voltage) is rectangular in order to obtain a triangular induction. The $B$ harmonic distortion  is apparent for high induction levels (near saturation).}
		\label{figure_deformB_sans_asserv}
\end{figure}

Otherwise, the measurements must be ensured in conformity with the international standards \cite{IEC}. Of the standard requirements to be considered, the induction waveform must ensure a sinusoidal waveform with respect to time (in the frequency range) with a form factor $FF < 1\%$. $FF$ given by (\ref{eq:fact_forme_pourcent}) the parameter quantifying the harmonic distorsion of the signal. This factor corresponds to the ratio between the RMS signal value over the average value of the rectified signal. $FF_{Theoretical}$ in (\ref{eq:fact_forme_pourcent}) is  equal to $\frac{\pi}{2\sqrt{2}} = 1, 1107$ for a purely sinusoidal signal, and equal to 1 for a purely square signal. We express it in \%  for commodity and we generalize it for any signal waveform (\ref{eq:fact_forme_pourcent}). Thus, $0\%$ of $FF$ means that the signal has the expected theoretical waveform. Keeping $FF < 1\%$ becomes very difficult especially at high induction levels given, on one hand, the high nonlinear nature of the materials and on the other hand, the Epstein frame itself which could be considered as a nonlinear transformer (i.e., not a pure inductance). For instance, Fig. \ref{figure_deformB_sans_asserv} shows an example of Epstein frame measurements made without feedback on a non-oriented SiFe steel under 5 Hz quasi-triangular waveform near saturation. The harmonic distorsion of $B$ (see the solid line curve) is clearly visible.

\begin{equation}\label{eq:fact_forme_pourcent}
FF(\%) = 100 \times \left( \frac{ \frac{v_{B_{eff}}}{<|v_{B}|>} - FF_{Theoretical}}{FF_{Theoretical}} \right)
\end{equation}
where: $v_{B_{eff}}$ is the RMS scalar value measured on the $v_B$ voltage, and $<|v_{B}|>$ is the mean scalar value measured on the rectifier $|v_{B}|$ voltage. $FF(\%)$ is the form factor wich quantifies the distorsion.\\
\noindent To control the induction waveform $B$ (i.e., $v_B$), it is necessary to resort to an analog or digital feedback. The idea consists in modifying the voltage waveform $v_{input}$ sent to the power supply. In our case, the digital control is privileged because:
\begin{enumerate}
\item It allows a good reproducibility;
\item It allows low self-oscillations contrary to the analog electronic circuit which fails given the oscillations in the feedback loop, especially for high induction levels;
\item It "manages" the measured noise in the measurement chain;
\item The control of the measured induction $B$ (i.e., $v_B$) is possible taking into account numerically the compensation of the air flux. 
\end{enumerate}


A digital controller calculates the correct waveform to be applied in order to obtain the desired induction in the material using a recursive procedure \cite{Bertotti}, \cite{Fiorillo}. Such procedure can be applied eventually in real-time correction.\\

Certain waveform digital control methods for magnetic characterization setups are reported in the literature. Some authors proposed a solution based on harmonic compensation \cite{Zhang} to correct the induction waveforms. The latter is limited to sinusoidal waveforms or needs a complex correction (intrinsically unstable) of the phases associated to the signal harmonics \cite{zurek}. In \cite{Matsubara}, another approach is investigated based on the identification of the equivalent electric circuit elements of the magnetic circuit to be characterized.\\

We propose in this paper an original method for digital control of waveforms in magnetic measurement setups. This method is based on the model-free controller recently proposed in 2008 for which, we propose a variant that removes the use of the estimation of the numerical derivative of y (i.e $v_B$). The ability of this digital control method to control the induction waveforms for peak induction up to near saturation is fully described and illustrated by both simulation and measurements. The main objective is to ensure low harmonic distorsion of the waveforms.

\section{Model-Free Based Digital Control}

\subsection{From the PID regulator to the Model-Free Controller}
The proposed "model-free"-based controller has been built from the approach originally developed by Fliess \& Join \cite{fliess2}. 
It makes {\it a priori} more robust the classical control law of PID type regarding the incertitudes and drift of the controlled model. 
The main advantages of the method can be summarized as follows:
\begin{enumerate}
\item Its simplicity to use and implement in comparison with e.g. the robust control;
\item Its ability to stabilize a wide range of nonlinear dynamical systems and to maintain very interesting dynamical performances in presence of strong perturbations.
\end{enumerate}

Such control law needs only few adjustments for a (relatively) good efficiency. 
To increase the performances of the method, the parameters of the control law should be adjusted depending on the knowledge of the controlled dynamical system.
However, we can experimentally verify that, for a given set of parameters, the dynamic and static performances remain very interesting in a wide range of (derived) parameters 
of the controlled system.\\

\begin{figure}[!htb]
	\centering
		\includegraphics[width=10.0cm]{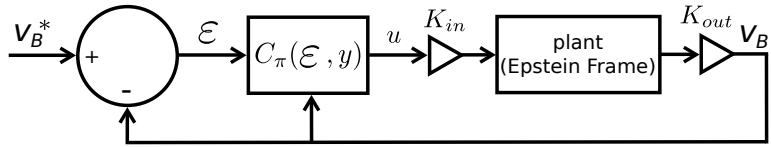}
		\caption{Proposed model-free scheme to control the Epstein Frame.}
		\label{fig:CSM_gen}
\end{figure}

Initially, the model-free control methodology has been successfully applied to many mechanical and electrical systems (see \cite{fliess2} for a non-exhaustive 
list of applications including power electronic converters \cite{Michel}. Described as a self-tuning controller in \cite{kumar}, this method has been improved
in \cite{Michel_arxiv} and last advances include the "extremum seeking control" for nonlinear systems.\\
The proposed improved model-free controller can be seen as a variant of the original discrete model-free control law: the computation of the numerical derivatives in 
the original approach \cite{fliess2} is substituted in the present proposed method by an initialization function that makes the controller more robust especially to the noise. 
The aim is to fit the response $v_{B}$ (in Fig. \ref{fig:CSM_gen}) with the output reference $v_{B}^{\star}$ by minimizing, for example, the square tracking error.

\begin{figure}[!htbp]
\centering
\subfloat[5 Hz]{\includegraphics[width=7cm]{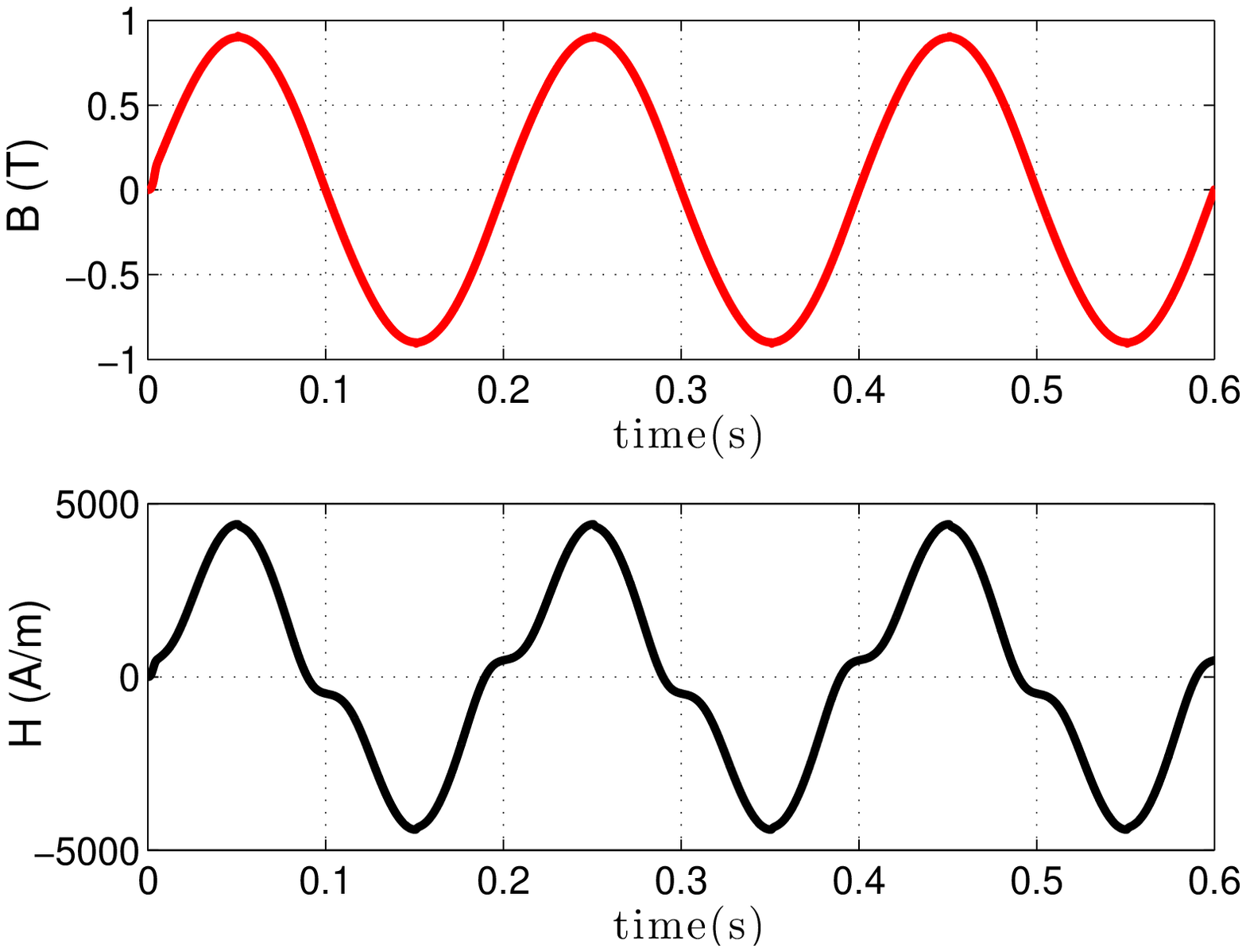}
\label{fig:fig_sim_10Hz}}
\subfloat[500 Hz]{\includegraphics[width=7cm]{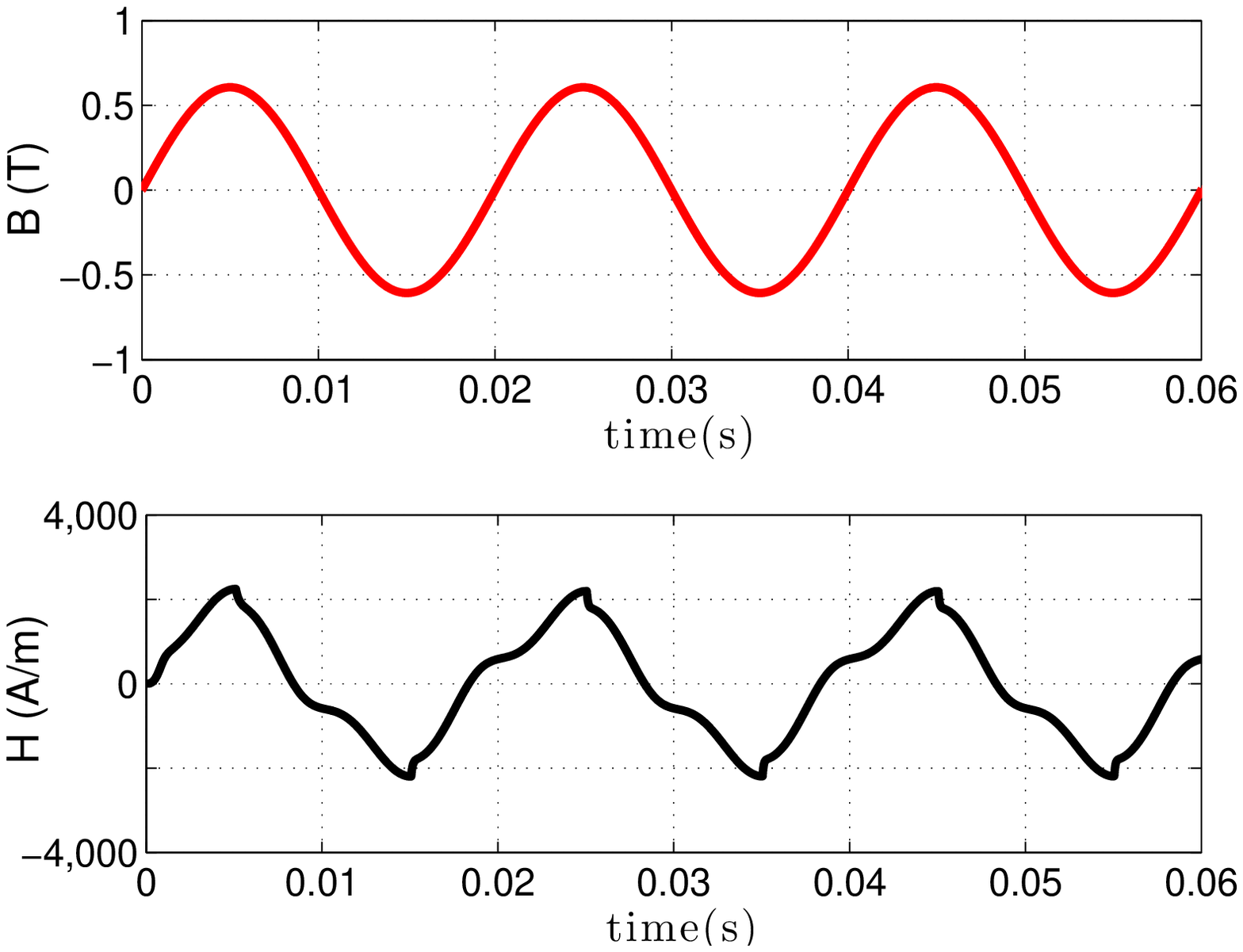}
\label{fig:fig_sim_500Hz}}
\caption{Examples of simulation of $u = H \, $[A/m] and $ v_B = B\, $[T] as a function of time in quasi-static (5 Hz) and dynamic (500 Hz) regimes.}
\label{fig:fig_sim}
\end{figure}
In simulation, an adjustment of  the controller parameters, depending on the model to be controlled, can be performed by using an optimization procedure. 
The aim is to fit the response $v_B$ (in Fig. \ref{fig:CSM_gen}) with the output reference ${v_B}^{\star}$ by minimizing, for example, the square tracking error. 
To do so, it is possible, for example, to use a meta-heuristic method, such as the standard simulated annealing \cite{Eglese} or a strategy based on 
derivative-free optimization such as, the brute force optimization (BFO) \cite{BFO} recently introduced.

\subsection{Technical Description}
Consider the control scheme depicted in Fig. \ref{fig:CSM_gen} where $\mathcal{C}_{\pi}$ is the proposed model-free $\&$ derivative-free controller 
(that is a modified version of the original Fliess controller).
For each discrete time $t_k, \, k \in \mathbb{N}^*$, the discrete controller $\mathcal{C}_{\pi}$ is defined by:
\begin{equation} \label{eq:iPI_discret_nm_eq}
u_k =  \left. \int_0^t K_i \varepsilon_{k-1} d \, \tau \right|_{k-1} \left\{ u_{k-1}^{i} + {K_p} ( k_{\alpha} e^{-k_{\beta}\, k} - v_{B_{k-1}}) \right\}
\end{equation}
\noindent
where: $u_k$ is the output of the controller that drives the input voltage of the Epstein frame; $v_{B}^{\ast}$ is the output reference trajectory; $K_p$ and $K_i$ are real 
positive tuning gains (subscript "p" and "i" denote proportional and integrative terms respectively); $\varepsilon_{k-1} = v_{{B}_{k-1}}^{\ast} - v_{B_{k-1}}$ is the tracking 
error; 
$k_{\alpha} e^{-k_{\beta} \, k}$ is an initialization function where $k_{\alpha}$ and $k_{\beta}$ are real constants; in practice, the integral part is discretized 
using e.g. Riemann sums. The internal recursion on $u_{k-1}^{i}$ is defined such as: $u_{k}^{i} = u_{k-1}^{i} + K_{p}(k_{\alpha} e^{-k_{\beta}\, k} - v_{B_{k-1}})$.

The current theoretical investigations aim to justify the nature of the initialization function (as a decreasing exponential function of the time) using applied variational calculus. 

\subsection{Practical Implementation}
As mentioned before, Fig. \ref{CaracBench} shows a general scheme of the Epstein frame setup used in magnetic measurements. Labview software is used to drive the setup. The purpose of the control law implementation is to control the input voltage $u$ such as the output voltage of the Epstein Frame $v_B$ remains "as close as possible" to a desired waveform reference. The controller $\mathcal{C}_{\pi}$ drives the linear voltage supply $u$, the measured voltage $v_B$ is recorded through an acquisition board to provide the feedback to the controller $\mathcal{C}_{\pi}$.

\section{Simulation}
In this section, we evaluate the dynamic performances of the proposed model-free controller in simulation. For simplicity purpose, 
the Epstein frame in Fig. \ref{CaracBench} is replaced by a nonlinear dynamic system $f_{BH}(i_H)$ describing the hysteresis loop $B(H)$ as a function of the 
excitation current $i_H$ such as in (\ref{eq:sim}). See \cite{Michel_arxiv} for more simulation results. 
\begin{equation}
\label{eq:sim}
v_B = B_{JA} (u).
\end{equation}
The function $B_{JA}$ considered in our case is the Jiles-Atherton (JA) model \cite{Jiles}, both static and dynamic versions \cite{Jiles1}, \cite{Jiles2}. Notice that the JA parameters were identified from a measured static major loop \cite{Messal}.\\
In  simulation,  the  system (\ref{eq:sim}) is  controlled  in  its  ”original  formulation”  without  any modification / linearization.

The BFO \cite{BFO} has been used to adjust the parameters of the $\mathcal{C}_{\pi}$ controller in such manner that the output response waveform $v_B$ remains "as close to" a sine waveform.\\
Fig. \ref{fig:fig_sim} depicts the input $u$ and the output $v_B$ of the control loop as a function of time in both quasi-static (Fig. \ref{fig:fig_sim_10Hz}) and dynamic (Fig. \ref{fig:fig_sim_500Hz}) regimes for a sinusoidal output reference ${v_B}^*$. 
The simulation results show that the closed-loop gives interesting dynamic performances, in particular, the form factor $FF$ is less than $1$\%.
\begin{figure}[!htbp]
\centering
\subfloat[10 Hz]{\includegraphics[width=7cm]{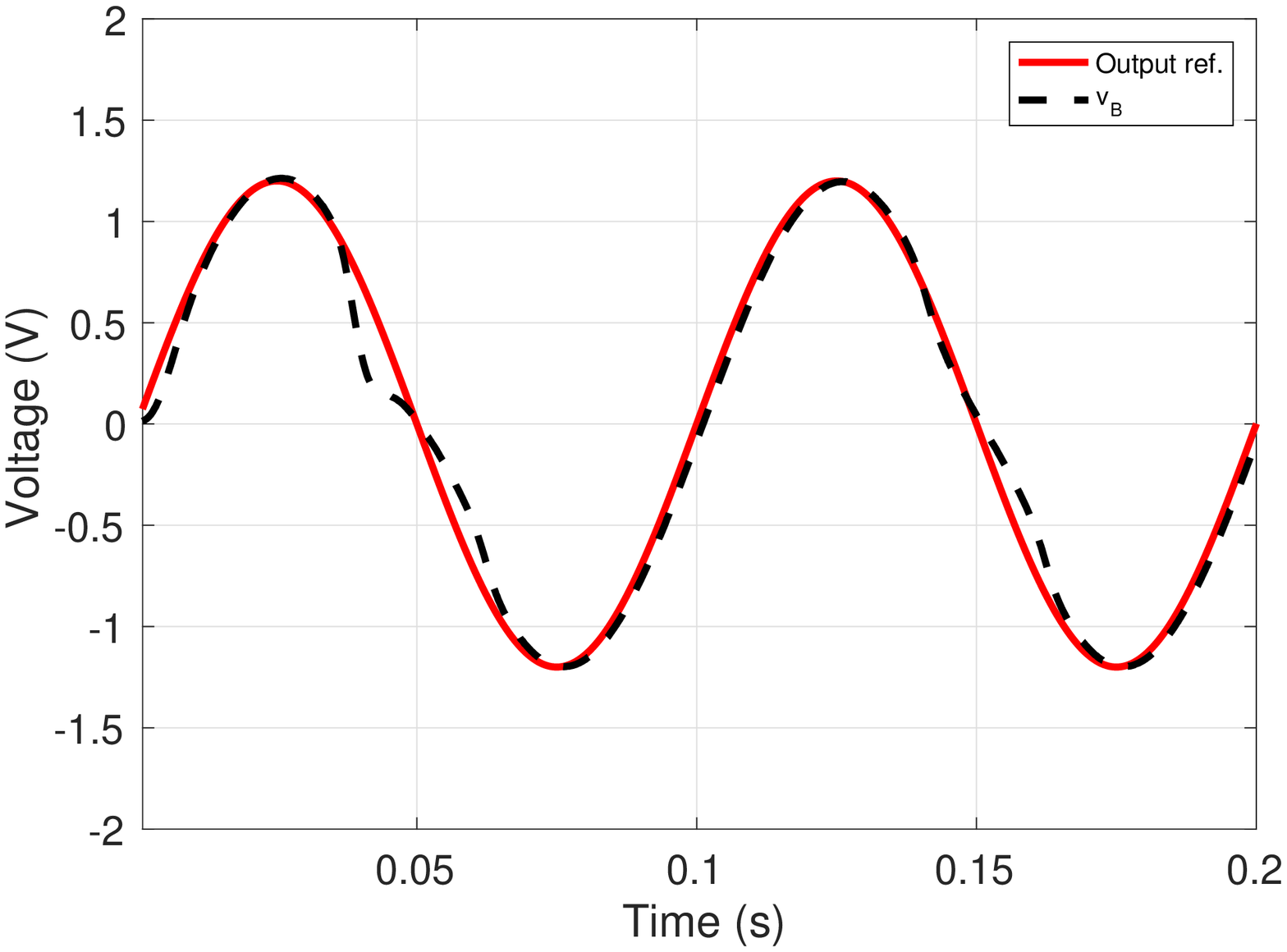}
\label{v_B_10Hz}}
\subfloat[50 Hz]{\includegraphics[width=7cm]{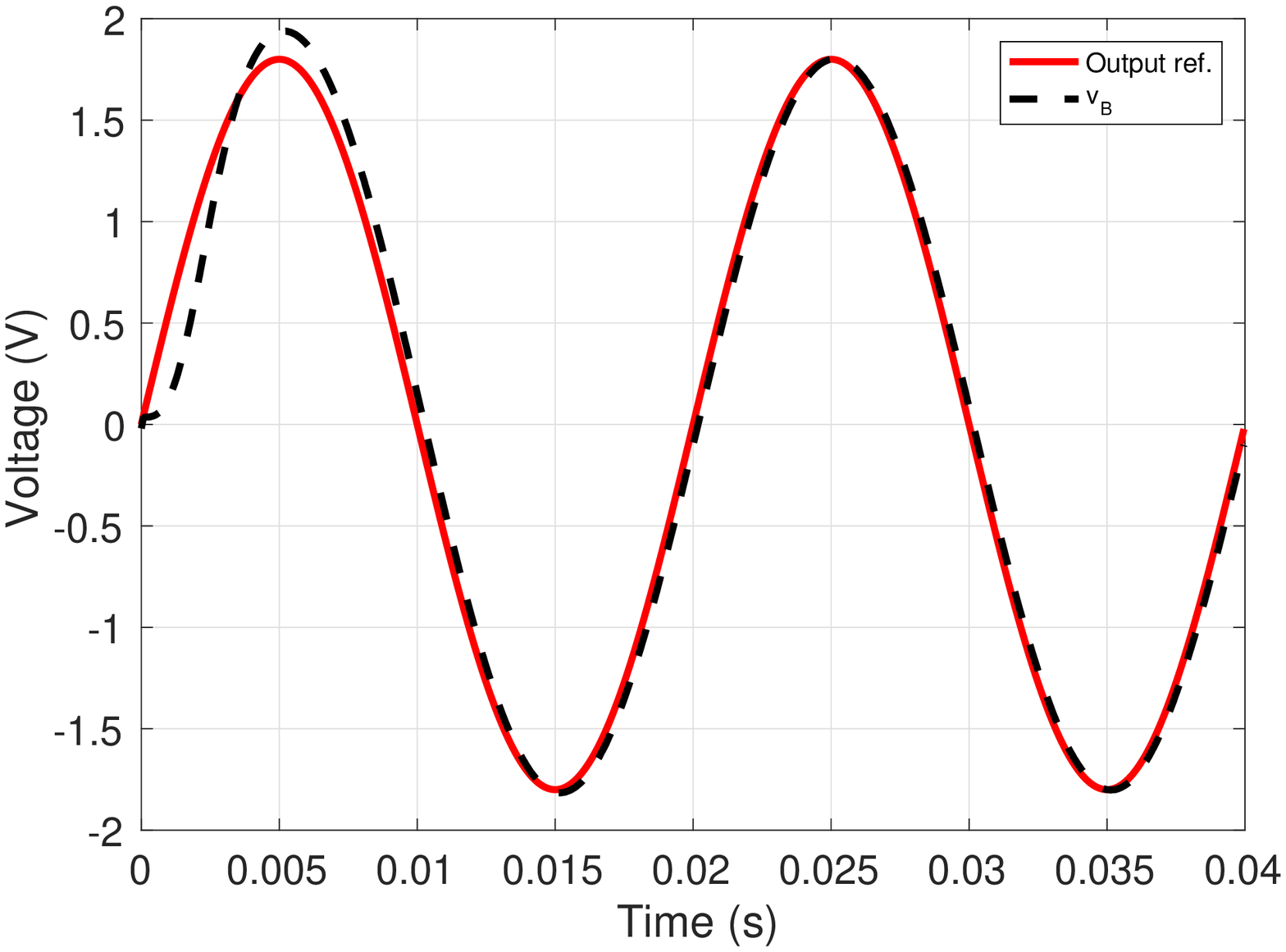}
\label{v_B_50Hz}}
\caption{Experimental performances of the proposed control in quasi-static (10 Hz) and dynamic (50 Hz) regimes.}
\label{Exp_10&50Hz}
\end{figure}

\section{Experimental Implementation}
\subsection{Sine Output References}
\begin{figure}[!htbp]
\centering
\includegraphics[width=8.7cm]{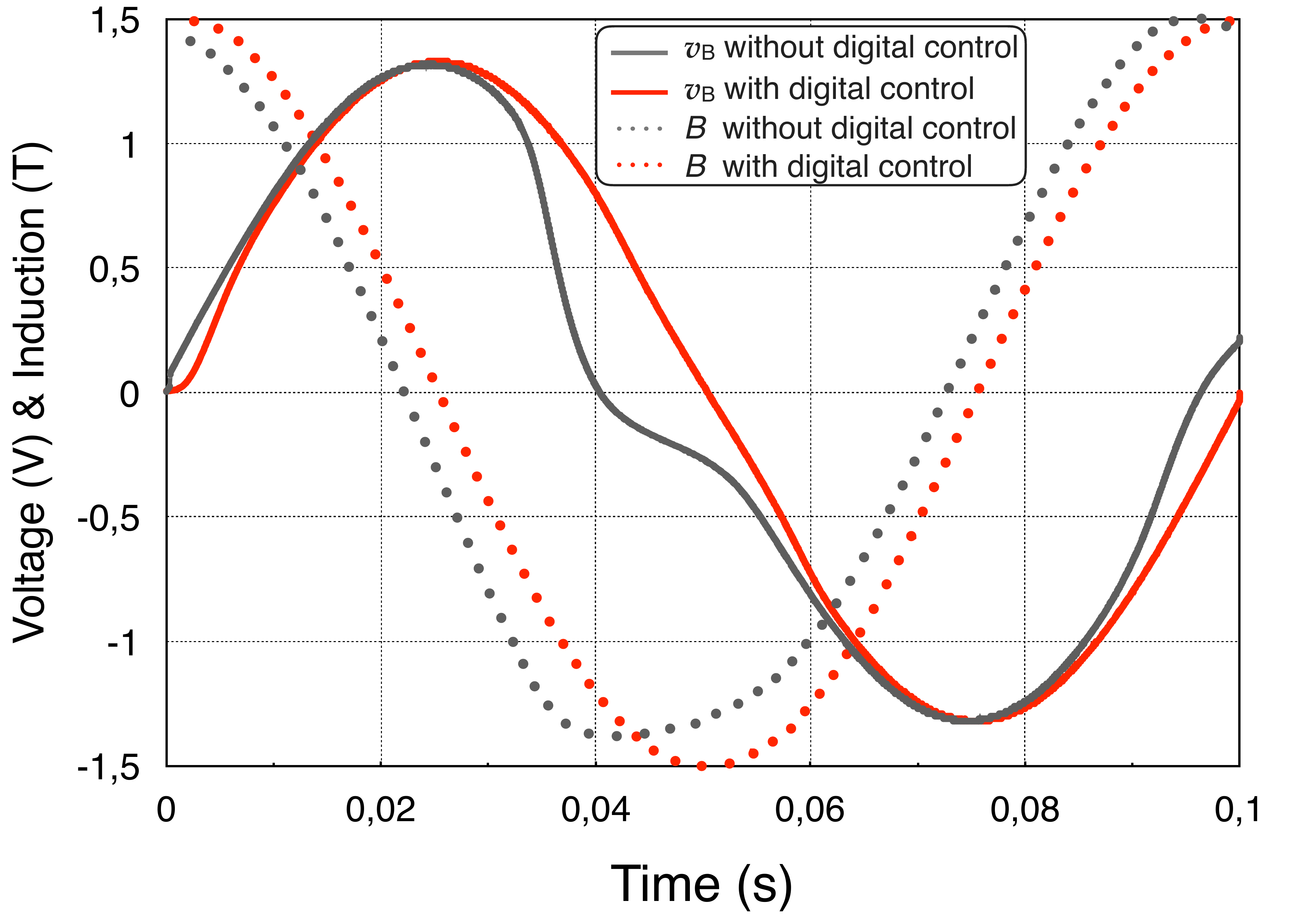}
\caption{Measurements performed with the proposed digital control at 10 Hz in comparison with that of no control. The form factor $FF(\%)$ given by (\ref{eq:fact_forme_pourcent}) but applied to $B(t)$ signal is clearly improved with the digital control (from 6.1\% without control to 0.53\% with control).}
\label{V2_etB_vs_t}
\end{figure}

In sine waveform case, we distinguish two main studies: 
\begin{enumerate}
\item The control is immediately started after the setup is powered;
\item The control is started after $n$ initialization cycles in order to damp the transient effects of the setup. 
\end{enumerate}
Both cases have been investigated and we focus on the results of the last one that gives better performances. Fig. \ref{Exp_10&50Hz} shows the performances of the proposed control on the magnetic characterization setup at 10 Hz (Fig. \ref{v_B_10Hz}) and 50 Hz (Fig. \ref{v_B_50Hz}). \\
The performances, in comparison with the experiments performed without control are shown in Fig. \ref{V2_etB_vs_t}. Notice that in open-loop (without control), a voltage drop occures (e.g., between $[0.03 - 0.05$ s$]$ for 10 Hz test). The latter is damped by the controller in closed-loop.

\subsection{Square Output References}
\begin{figure}[!htbp]
\centering
\subfloat[Controlled $v_B $ and $B$]{\includegraphics[width=7cm]{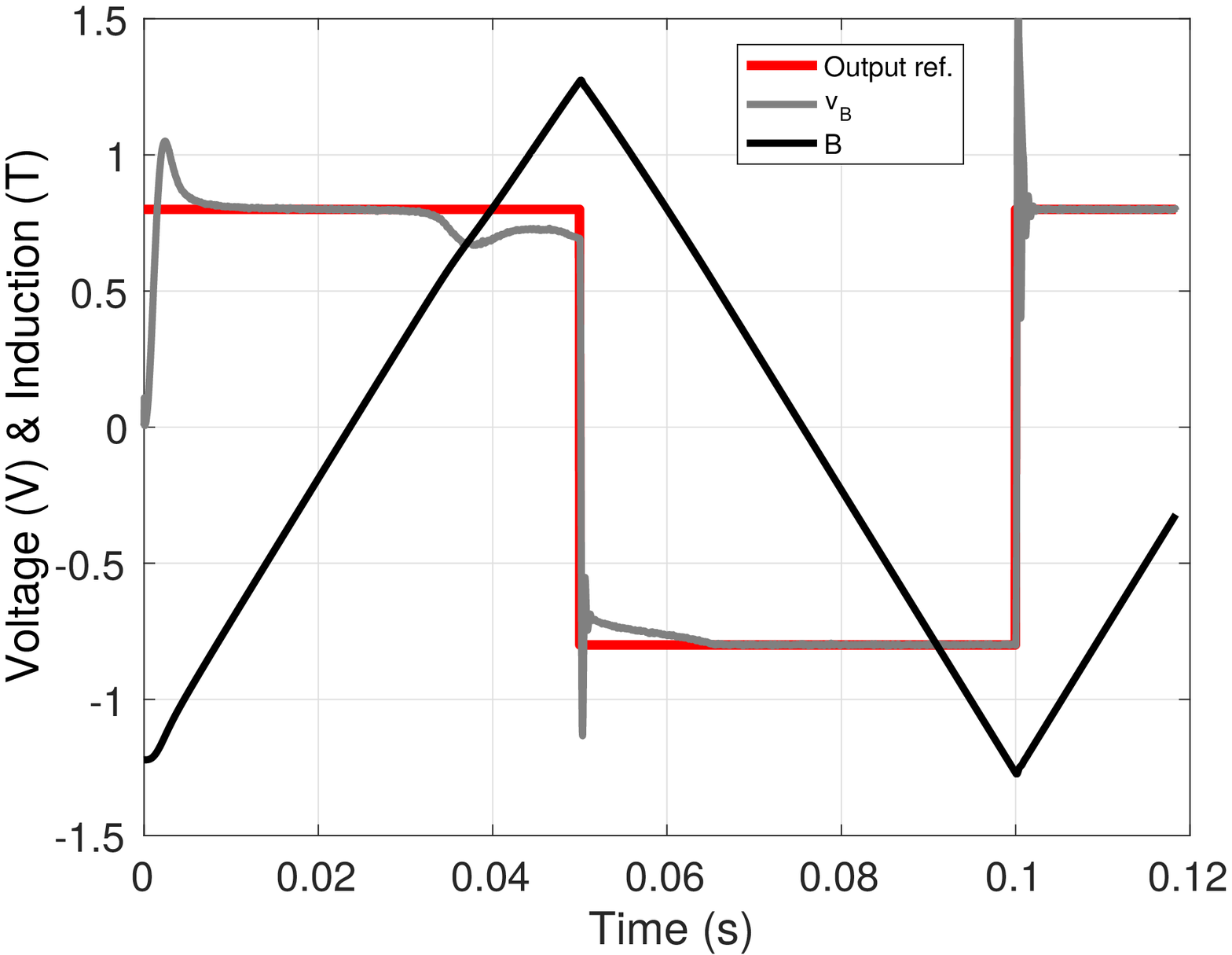}
\label{TriangleControl_v_B_10Hz}}
\subfloat[$u$]{\includegraphics[width=7cm]{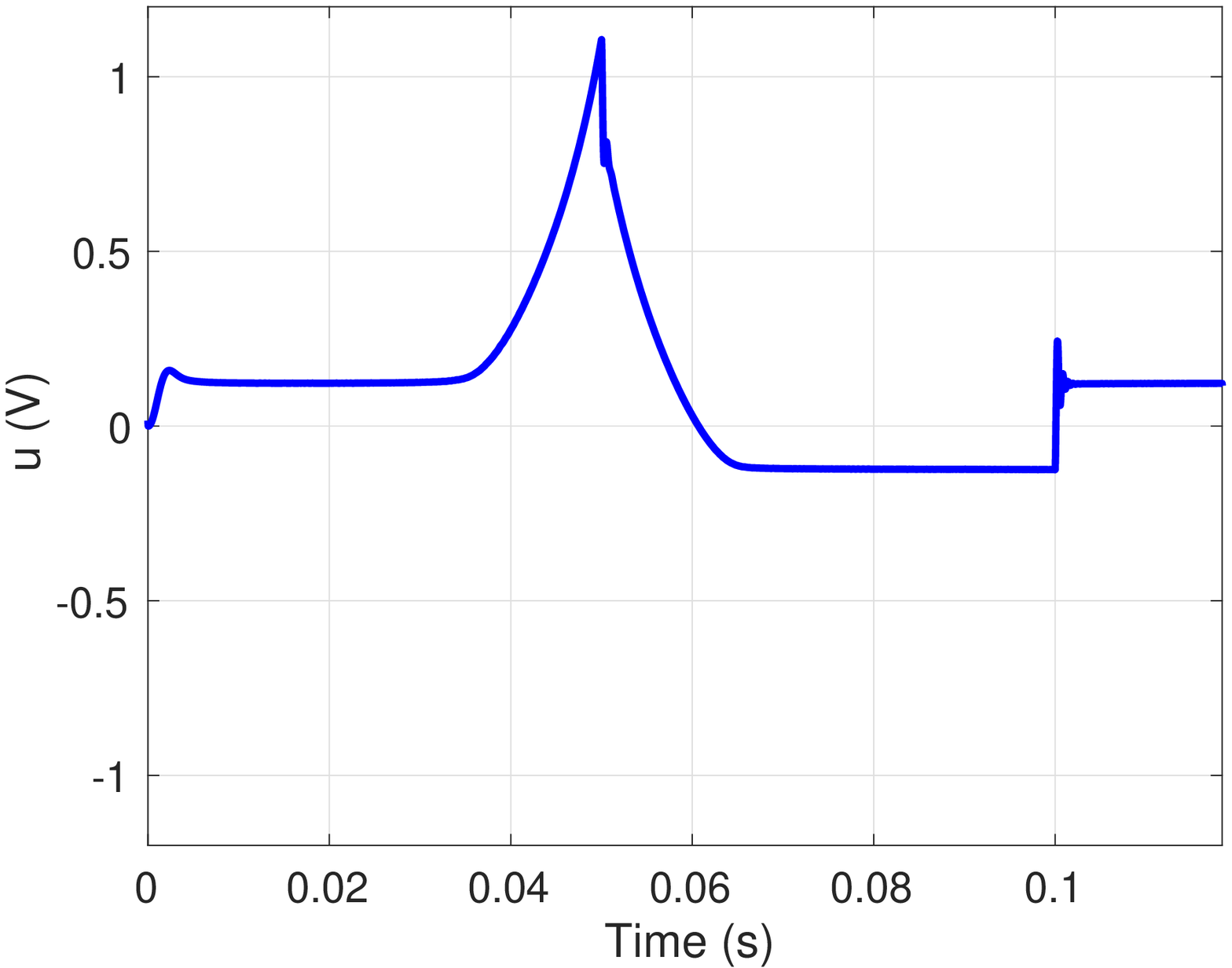}
\label{u_triangle_10Hz}}
\caption{Measurements performed with the proposed digital control at 10 Hz for a square output reference.}
\label{TriangleControl_10Hz}
\end{figure}

\begin{figure}[!htbp]
\centering
\subfloat[Controlled $v_B $ and $B$]{\includegraphics[width=7cm]{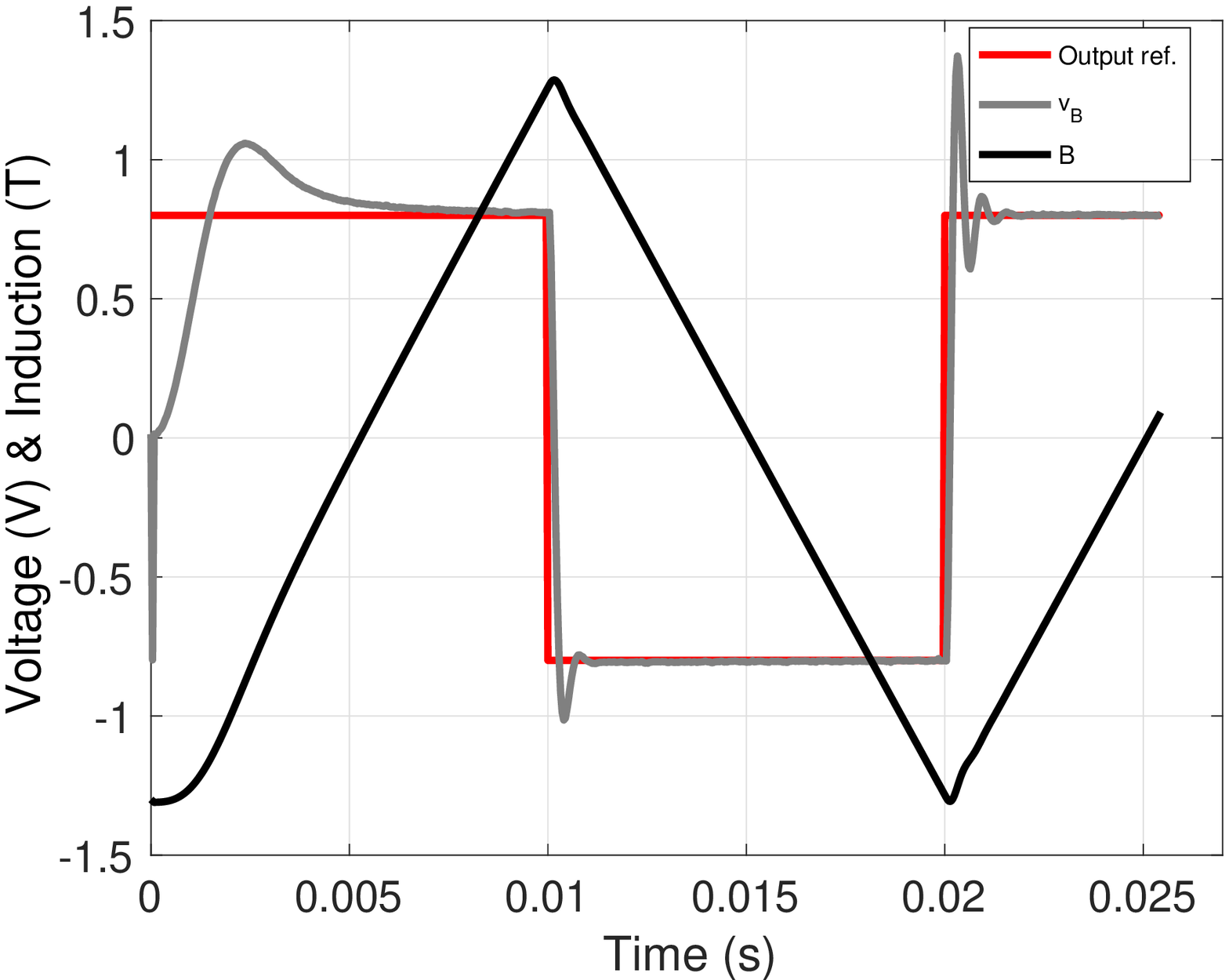}
\label{TriangleControl_v_B_50Hz}}
\subfloat[$u$]{\includegraphics[width=7cm]{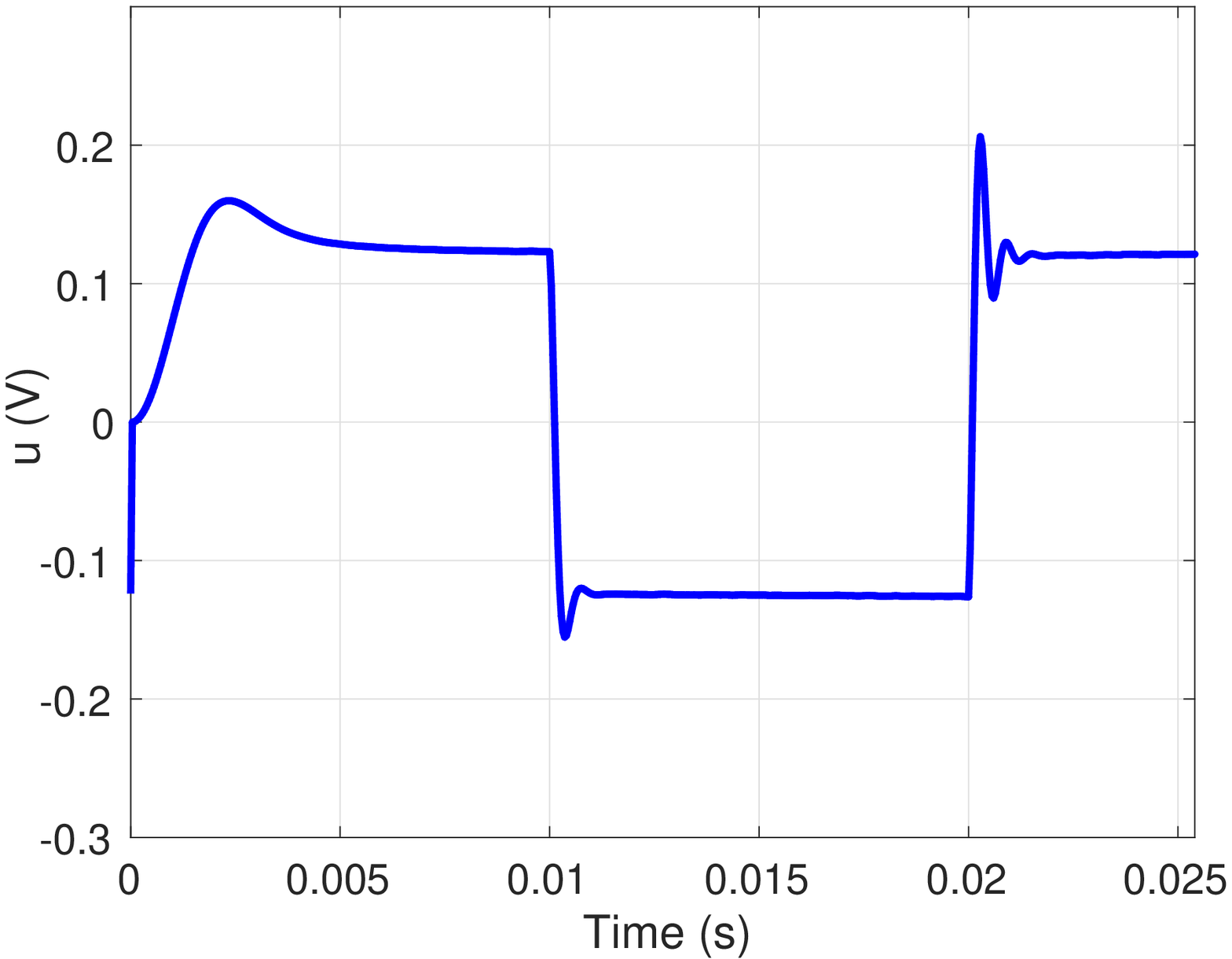}
\label{u_triangle}}
\caption{Measurements performed with the proposed digital control at 50 Hz for a square output reference.}
\label{TriangleControl_50Hz}
\end{figure}

Fig. \ref{TriangleControl_10Hz} and Fig. \ref{TriangleControl_50Hz} show the technical performances of the proposed digital control in square output reference case respectively measured at 10 and 50 Hz. One observes the improvement brought by the proposed digital control in comparison with the measurements made in open-loop given in Fig. \ref{figure_deformB_sans_asserv}. Indeed, the measured induction $B$ (obtained after integrating the secondary voltage  $v_B$) presents a quasi-triangular waveform.\\

As illustrated, for example in Fig. \ref{u_triangle_10Hz} (i.e., at 10 Hz), we observe that $u$ strongly increases before the end of the half period. This is due to $v_B$ voltage drop occuring at the same time. Like in sine waveform case, the controller compensates the $v_B$ voltage drop.\\
\subsection{Symmetrization}
\begin{figure}[!htpb]
\centering
\subfloat[Measured $i_H$ and controlled $v_B $]{\includegraphics[width=7cm]{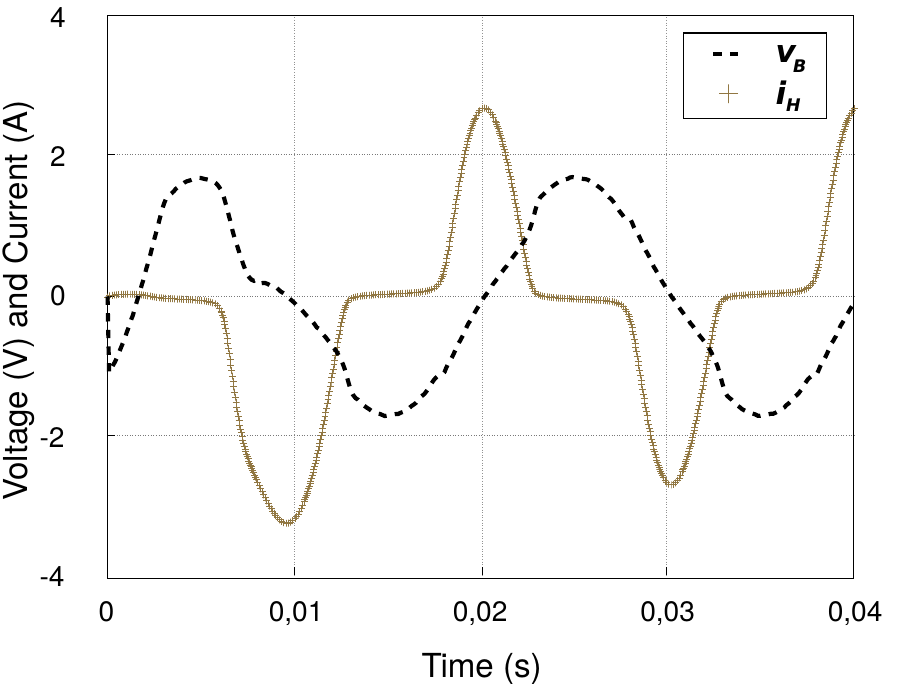}
\label{Sym1}}
\subfloat[Measured hysteresis loop]{\includegraphics[width=7cm]{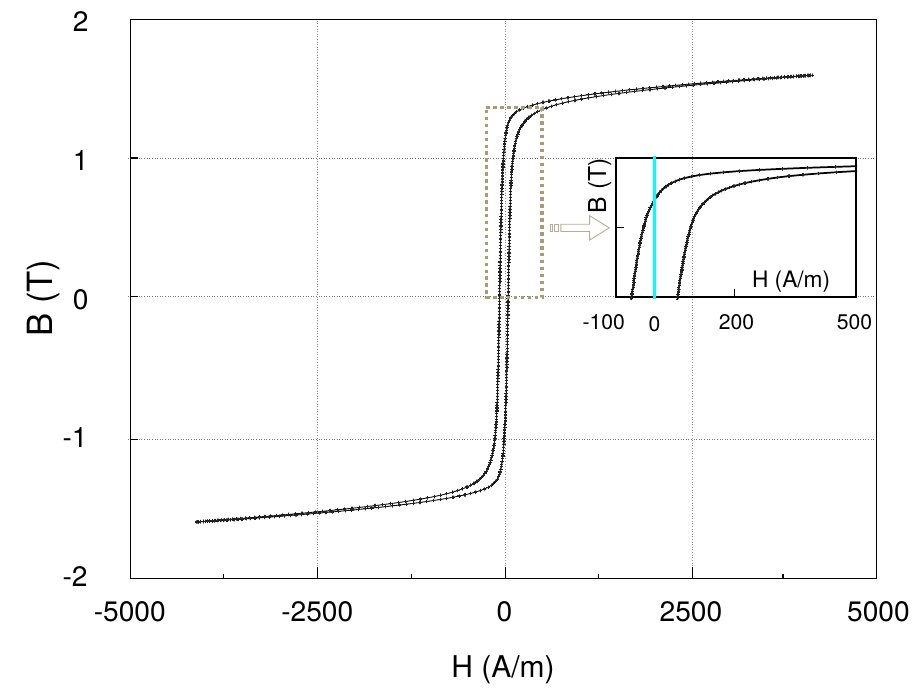}
\label{Sym2}}
\caption{Measurements performed with the proposed digital control at 50 Hz with symmetrization.}
\label{Sym}
\end{figure}

The measurments presented above do not preserve the symmetry of the current $i_H$. This is unsuitable for magnetic meausrements because it introduces a power dc component which magnetize the material. To overcome this limitation, we introduced a symmetrization technique which aims to adjust the voltage $u$ in addition to the modifications induced by the controller.
For example, Fig. \ref{Sym1} shows the measured performances of the proposed control -with symmetrization- in sine output reference case at 50 Hz. We can observe the symmetry of the current $i_H$ and thus the applied excitation field $H$ in Fig. \ref{Sym2} displaying the measured hysteresis loop and zoomed-in view of the zone around $H_c$. One observes that the resulting measured hysteresis loop does not present any offset.
Even if the form factor $FF(\%)$ given by (\ref{eq:fact_forme_pourcent}) is increased in comparison with the case without symmetrization, the measurements still in accordance with standards with a good form factor (e.g., 0.7\% for the given example in Fig. \ref{Sym1}). 
\section{Conclusion and Prospects}
In this paper, we have presented an original implementation of a digital control for a characterization bench of soft magnetic materials based on the model-free approach. The results showed good dynamic performances regarding both sinusoidal and square output references. The form factor $FF(\%)$, which quantifies the harmonic distorsion of $v_B$, is greatly improved with the proposed digital control (at least by 10).  The measured induction $B$ (obtained after integrating the $v_B$ voltage) presents respectively a quasi-sine and triangular waveform.\\

Future developments should focus on the adaptive adjustment of the control coefficients as a  function of the considered material; the excitation frequency and  the output reference. In addition, a recent work which is currently in progress will compare the model-free based digital control with another waveform control method based on the knowledge of the characterization setup through an electric equivalent circuit model of the Epstein frame. The comparison will focus on the performances in terms of rapidity and accuracy in different working conditions.
Furthermore, other materials (such as amorphous or nanocrystalline materials, etc) with different geometries (such as torus, single sheet tester) can be tested. The proposed digital control can also be applied for rotating magnetic field experiments or biaxial characterizations.

\section*{Acknowledgment}
This work was supported by Renault-SAS, Guyancourt, France, through the Program COCTEL financed by the French Agency for Environment and Energy Management ADEME.


\end{document}